# A few equalities involving integrals of the logarithm of the Riemann $\varsigma$-function and equivalent to the Riemann hypothesis

Sergey K. Sekatskii, Stefano Beltraminelli, and Danilo Merlini

**Abstract.** Using a generalized Littlewood theorem concerning integrals of the logarithm of analytical functions, we have established a few equalities involving integrals of the logarithm of the Riemann $\varsigma$-function and have rigorously proven that they are equivalent to the Riemann hypothesis. Separate consideration for imaginary and real parts of these equalities, which deal correspondingly with the integrals of the logarithm of the module of the Riemann function and with the integrals of its argument is given. Preliminary results of the numerical research performed using these equalities to test the Riemann hypothesis are presented.

## 1. Introduction

Search for criteria equivalent to the Riemann hypothesis (RH) and involving integrals of the logarithm of the Riemann $\zeta$-function (see e.g. [1] for standard definitions and discussion of the general properties of this function) started with the paper of Wang [2] and then was continued by Volchkov [3], Balazard – Saias - Yor [4] and later on by Merlini [5] and Merlini *et al*. [6]. In this Note we would like to develop further and to generalize the same line of research which has been started in [5, 6], and to introduce a few equalities involving the integrals of *ln($\zeta$)* and equivalent to the RH.

In essence, our approach is a slight generalization of the Littlewood theorem concerning contour integrals of the logarithm of analytical function (see, for example [7]), and an application of this generalization to *ln($\zeta$)*. Hence the present note is to a large extent complementary to the earlier papers [5, 6].

## 2. Generalization of the Littlewood theorem concerning integrals of the logarithm of an analytical function

Throughout the paper we use the notation $z = x + iy$ and/or $z = \sigma + it$ on the equal footing. Let us consider functions *f(z), g(z)* and let *F(z)=ln(f(z))*. Let us introduce and analyze an integral $\int_C F(z)g(z)dz$ along the contour *C* which is the rectangle formed by the broken line composed by four straight line segments connecting the vertices with the coordinates $X_1 + iY_1$, $X_1 + iY_2$, $X_2 + iY_2$, $X_2 + iY_1$; we suppose that $X_1 < X_2$, $Y_1 < Y_2$ and that *f(z)* is analytic and non-zero on *C* and meromorphic inside it, and that *g(z)* is analytic on *C* and meromorphic inside it. Assume first that the poles and zeroes of the functions *F(z), g(z)* do not coincide and that the function *f(z)* has no poles and only one simple zero located at the point $X+iY$ in the interior of the contour. Then let us do the cut along the straight line segment $X_1+iY$, $X+iY$ and consider a new contour $C'$ which is the initial contour $C$ plus the cut indenting the point $X+iY$ (see Fig. 1; cf. the proof of the Littlewood theorem in [7]). Cauchy residue theorem can be applied to the contour $C'$ which means that the value of $\int_{C'} F(z)g(z)dz$ is $2\pi i \sum_{\rho} F(\rho)res(g(\rho)))$, where the sum is over all simple poles $\rho$ of the function $g(z)$ lying inside the contour.

An appropriate choice of the branches of the logarithm function assures that the difference between two branches of the logarithm function appearing after the integration path indents the point $X+iY$ is $2\pi i$ (see below for an exact definition), which means that the value of an integral along the initial contour $C$ is $2\pi i (\sum_{\rho} F(\rho)res(g(\rho))) - \int_{X_1+iY}^{X+iY} g(z)dz)$ (again, compare with the proof of the Littlewood theorem in [7]. Remind, that the integral is taken here along the straight segment parallel to the real axis hence for this integral $Im(z)=const=Y$ and $dz=dx$). Analogously, a single pole of $f(z)$ function occurring at (other) point $X+iY$ contributes $2\pi i \int_{X_1+iY}^{X+iY} g(z)dz$ to the contour integral value. Thus by summing all terms arising from different zeroes and poles of the function $f(z)$ one obtains the following generalization of the Littlewood theorem:

THEOREM 1. Let $C$ denotes the rectangle bounded by the lines $x=X_1$, $x=X_2$, $y=Y_1$, $y=Y_2$ where $X_1<X_2$, $Y_1<Y_2$ and let $f(z)$ be analytic and non-zero on $C$ and meromorphic inside it, let also $g(z)$ is analytic on $C$ and meromorphic inside it. Let $F(z)=ln(f(z))$, the logarithm being defined as follows: we start with a particular determination on $x = X_2$, and obtain the value at other points by continuous variation along $y=const$ from $\ln(X_2 + iy)$. If, however, this path would cross a zero or pole of *f(z)*, we take *F(z)* to be $F(z \pm i0)$ according as we approach the path from above or below. Let also the poles and zeroes of the functions *f(z), g(z)* do not coincide.

$$\text{Then } \int_C F(z)g(z)dz = 2\pi i (\sum_{\rho_g} res(g(\rho_g))F(\rho_g)) - \sum_{\rho_f^0} \int_{X_1+iY_\rho^0}^{X_\rho^0+iY_\rho^0} g(z)dz + \sum_{\rho_f^{pol}} \int_{X_1+iY_\rho^{pol}}^{X_\rho^{pol}+iY_\rho^{pol}} g(z)dz)$$

where the sum is over all $\rho_g$ which are simple poles of the function *g(z)* lying inside *C*, all $\rho_f^0 = X_\rho^0 + iY_\rho^0$ which are zeroes of the function *f(z)* counted taking into account their multiplicities (that is the corresponding term is multiplied by *m* for a zero of the order *m*) and which lye inside *C*, and all $\rho_f^{pol} = X_\rho^{pol} + iY_\rho^{pol}$ which are poles of the function *f(z)* counted taking into account their multiplicities and which lye inside *C*. For this is true all relevant integrals in the right hand side of the equality should exist.

Of course, exact nature of the contour is actually irrelevant for the theorem as it is for the Littlewood theorem. As for the last remark concerning the existence of integrals taken along the segments $[X_1 + iY_\rho, X_\rho + iY_\rho]$ parallel to the real axis (the inexistence might happen e.g. if a pole of the *g(z)* function occurs on the segment), the corresponding integration path often can be modified in such a manner that one will be able to handle these integrals. The case of the coincidence of poles and zeroes of the functions *f(z), g(z)* often does not pose real problems and can be easily considered. We will deal with a few such cases below.

**3. Applications of above theorem to the Riemann $\zeta$-function**

Now we are in a position to put $f(z)= \zeta(z)$, $g(z) = \dfrac{1}{a^2 - (z-b)^2}$ where for a moment $a \neq 0$, $b$ are arbitrary complex numbers. We consider the contour $C$ lying in the half-plane $Re(z)>0$ (indeed, this condition is not important but we want to avoid all unnecessary complications related with the trivial Riemann zeroes and asymptotic of the Riemann function for $Re(z)<0$) and shall search for the value of an integral $I = \displaystyle\int_C \dfrac{\ln(\zeta(z))}{a^2 - (z-b)^2} dz$. Function $g(z)$ has single poles at the points $z_1 = b+a$, $z_2 = b-a$, and these poles contribute correspondingly $-\dfrac{i\pi \ln(\zeta(b+a))}{a}$, $\dfrac{i\pi \ln(\zeta(b-a))}{a}$ to the contour integral provided, of course, that these same poles lie inside the contour. Using an identity $\displaystyle\int \dfrac{dz}{a^2 - (z-b)^2} = \dfrac{1}{2a} \ln\left(\dfrac{a+z-b}{a-(z-b)}\right)$, one obtains the following contribution to the contour integral arising from any Riemann function zero $\rho_k = \sigma_k + it_k$ lying inside the contour:

$$I_k = -2\pi i \int_{X_1+it_k}^{\sigma_k+it_k} \dfrac{dz}{a^2 - (z-b)^2} = -\dfrac{\pi i}{a}\left(\ln\left(\dfrac{a+\sigma_k+it_k-b}{a-\sigma_k-it_k+b}\right) - \ln\left(\dfrac{a+X_1+it_k-b}{a-X_1-it_k+b}\right)\right) \quad (1)$$

To find the value of $I$, we need to sum up such contributions $I_k$ for all zeroes of the Riemann function lying inside the contour taking into account their multiplicities; that is one needs to multiply the contribution (1) by $m$ for a zero of the order $m$. Finally, the simple pole $z=1$ of the Riemann $\zeta$-function, if it lies inside the contour, adds

$$I_{pol} = \dfrac{\pi i}{a}\left(\ln\left(\dfrac{a+1-b}{a-1+b}\right) - \ln\left(\dfrac{a+X_1-b}{a-X_1+b}\right)\right)$$ to the integral value.

The above expressions enable to calculate the value of any integral of the type $I = \int_C \frac{\ln(\varsigma(z))}{a^2 - (z-b)^2} dz$, and below we use them for the Riemann function research. Now let us consider a few particular interesting cases.

Let $g(z) = \frac{1}{a^2 - (z-b)^2}$ with $a$, $b$ real positive and $1 > b > 1/2$, and let the contour $C$ is formed by the vertices $b - iY$, $b + iY$, $b + X + iY$, $b + X - iY$; let also $X, Y \to +\infty$. Known asymptotic properties of the Riemann function for large argument values assure that the value of the contour integral taken along the "external" three straight lines containing at least one of the points $b + X + iY$, $b + X - iY$ tends to zero when $X, Y \to \infty$ (this has been shown for a particular case of $g(z)$ function, viz. $a=0$, $b=1/2$, already by Wang [2], and his consideration can be one-to-one repeated here for a more general case), and thus we obtain that the value of a contour integral tends to $I = -i \int_{-\infty}^{\infty} \frac{\ln(\varsigma(b+it))}{a^2 + t^2} dt$ (minus sign comes from the standard describing of the contour in the counterclockwise direction; integral is taken along the line $z=b+it$ where $dz = idt$ and $g(z) = \frac{1}{a^2 + t^2}$). According to our Theorem 1, the value of this integral is $\frac{i\pi}{a}\left(-\ln(\varsigma(a+b)) + \ln\left(\frac{a-b+1}{a+b-1}\right)\right)$ (the first term comes from the simple pole of $g(z)$ occurring at $z_1=a+b$ and thus definitely lying inside the contour while the second pole $z_2=b-a$ definitely lies outside it, and the second term comes from the simple pole of the Riemann function at $z=1$) plus the sum of terms involving all (possible) zeroes of the Riemann function occurring to the right of the line $Re(z)>b$. If RH holds true, there are no such zeroes due to our choice $b>1/2$ and obviously this last sum equals to zero. If there are such zeroes, their contribution equals to the sum of the terms $-\dfrac{\pi i}{a}\left(\ln\left(\dfrac{a+\sigma_k+it_k-b}{a-\sigma_k-it_k+b}\right)-\ln\left(\dfrac{a+it_k}{a-it_k}\right)\right)$ taken into account with the multiplicities of zeroes for all $\sigma_k > b$.

*3.1 Integrals involving logarithm of the module of the Riemann function*

Now it is instructive to consider separately real and imaginary parts of this equality. Due to the known symmetry properties of the Riemann function we have $\varsigma(z)=\varsigma^*(z^*)$; in particular, if $\rho_k=\sigma_k+it_k$ is a zero, $\rho_k=\sigma_k-it_k$ also is. This means that the above sum over the Riemann function zeroes has a purely imaginary value and the integral $I = -i\int_{-\infty}^{\infty}\dfrac{\ln(\varsigma(b+it))}{a^2+t^2}dt$ also has a purely imaginary value. This also means that

$\sum_{\sigma_k>b}\ln\left(\dfrac{a+it_k}{a-it_k}\right)=0$ and thus by division by $i$ we obtain:

$$\dfrac{a}{\pi}\int_{-\infty}^{\infty}\dfrac{\ln(|\varsigma(b+it)|)}{a^2+t^2}dt = \ln\left|\dfrac{\varsigma(a+b)\cdot(a+b-1)}{a-b+1}\right| + \sum_{\sigma_k>b}\ln\left|\dfrac{a+\sigma_k-b+it_k}{a-\sigma_k+b-it_k}\right| \quad (2)$$

This equality can be easily analyzed using the final expression $\sum_{\sigma_k>b}\ln\left|\dfrac{a+\sigma_k-b+it_k}{a-\sigma_k+b-it_k}\right|$ for the sum over the Riemann zeroes (which is written in (2)), but having in mind future applications let us use the initial form $\sum_{\sigma_k>b}\int_{b+it_k}^{\sigma_k+it_k}\dfrac{dz}{a^2-(z-b)^2}$. Variable change $z=p+b+it_k$ transforms each term of this sum into $\int_0^{\sigma_k-b}\dfrac{dp}{a^2-(p+it_k)^2}$. Its real part, which is currently under the consideration, is equal to $\int_0^{\sigma_k-b}\dfrac{a^2-p^2+t_k^2}{(a^2-p^2+t_k^2)^2+4p^2t_k^2}dp$ while its imaginary part

is equal to $\int_0^{\sigma_k - b} \frac{2ipt_k}{(a^2 - p^2 + t_k^2)^2 + 4p^2 t_k^2} dp$. Because this is well known that possible values of $t_k$ for the Riemann function zeroes with $\sigma_k > 1/2$, if they exist, are very large while $0 < p < 1/2$, the expression $a^2 - p^2 + t_k^2$ is always positive and hence all terms in the sum in the r.h.s. of (2) are also positive. Correspondingly, this sum can be equal to zero only if there are no Riemann function zeroes with $\sigma_k > b$. Evidently, zeroes with $\sigma_k = b$ contributes nothing to the sum and our approach is extended to include this case. Hence we have rigorously proven the following theorem earlier established by Merlini [5].

THEOREM 2. An equality $\frac{a}{\pi} \int_{-\infty}^{\infty} \frac{\ln(|\varsigma(b+it)|)}{a^2 + t^2} dt = \ln \left| \left( \frac{\varsigma(a+b) \cdot (a+b-1)}{a-b+1} \right) \right|$, where $a, b$ are arbitrary real positive numbers and $1 > b \geq 1/2$, holds true for some $b$ if and only if there are no Riemann function zeroes with $\sigma > b$. For $b=1/2$ this equality is equivalent to the Riemann hypothesis.

Taking $b>1$ and repeating the same steps as above we have unconditionally (because Riemann function has no zeroes when $\sigma>1$ and the point $z=1$ does not lie any more inside the integration contour) $\frac{a}{\pi} \int_{-\infty}^{\infty} \frac{\ln(|\varsigma(b+it)|)}{a^2 + t^2} dt = \ln(\varsigma(a+b))$.

Evident modification of condition (2) consists in the choice of $f(z) = \ln(\varsigma(z)(z-1))$ thus removing the pole at $z=1$. Repeating the same as above, one has the following

THEOREM 3. An equality

$$\frac{a}{\pi}\int_{-\infty}^{\infty}\frac{\ln(|\varsigma(b+it)\cdot(b-1+it)|)}{a^2+t^2}dt = \ln(|\varsigma(a+b)\cdot(a+b-1)|),$$ where *a, b* are arbitrary real positive numbers, $1 > b \geq 1/2$ and $a+b \neq 1$, holds true for some *b* if and only if there are no Riemann function zeroes with $\sigma > b$. For *b=1/2* this equality is equivalent to the Riemann hypothesis.

REMARK 1. Of course, one has simply $$\int_{-\infty}^{\infty}\frac{\ln|b-1+it|}{a^2+t^2}dt = \frac{\pi}{a}\ln(|a-b+1|),$$ see example #4.295.7 from Gradshtein and Ryzhik book [8] or apply directly the Theorem 1 from this Note, hence the theorem 3 can be seen as a trivial corollary of the Theorem 2.

Interesting simplification of condition (2) can be obtained by setting additionally *a+b=1*. It is easy to see that such a choice combines the simple pole of the function *g(z)* with the simple zero of the function $f(z) = \ln(\varsigma(z)(z-1))$ at *z=1* thus killing all peculiarities of the function under an integral sign, and one immediately obtains that the only contributions to the contour integral now arise from the Riemann zeroes lying to the right of the line *Re(z)>b*. Hence we have

THEOREM 4. An equality $$\int_{-\infty}^{\infty}\frac{\ln(|\varsigma(1/2+\alpha+it)(-1/2+\alpha+it)|)}{(1/2-\alpha)^2+t^2}dt = 0,$$ where $\alpha$ is an arbitrary real number $0 \leq \alpha < 1/2$, holds true for some $\alpha$ if and only if there are no Riemann function zeroes with $\sigma > 1/2+\alpha$. For $\alpha = 0$ this equality is equivalent to the Riemann hypothesis.

Of course, Theorem 4 follows from Theorem 3 when $a+b \to 1$.

## 3.2 Integrals involving the argument of the Riemann function

Additional equalities involving the integrals of the logarithm of the Riemann function can be established if we apply the Theorem 1 to the contour $C$ formed by the vertices $b$, $b+iY$, $b+X+iY$, $b+X$ with $1>b>1/2$, $X,Y \to +\infty$. For simplicity, we will consider only an integral $\int_C \frac{\ln(\varsigma(z) \cdot (z-1))}{(1/2-\alpha)^2 - (z-1/2-\alpha)^2} dz$; again, more general case can be easily investigated along the same lines. As before, the integrals along two "external" lines containing the point $b+X+iY$ tend to zero when $X,Y \to +\infty$. The integral taken along the real line is purely real, and hence the consideration of the imaginary part of the equality following from the Theorem 1 gives us nothing new. At the same time, the consideration of the real part of this equality gives:

$$\int_0^\infty \frac{\arg(\varsigma(1/2+\alpha+it)(-1/2+\alpha+it))}{(1/2-\alpha)^2+t^2} dt + \int_{1/2+\alpha}^\infty \frac{\ln(\varsigma(\sigma)(\sigma-1))}{(1-\sigma)(\sigma-2\alpha)} d\sigma =$$
$$2\pi \sum_{t_k>0, \sigma_k>1/2+\alpha} \int_0^{\sigma_k - 1/2 - \alpha} \frac{2pt_k}{((1/2-\alpha)^2 - p^2 + t_k^2)^2 + 4p^2 t_k^2} dp \quad (3)$$

To obtain this we used $\ln z = \ln|z| + i\arg(z)$ and expression $\int_0^{\sigma_k - b} \frac{2ipt_k}{(a^2 - p^2 + t_k^2)^2 + 4p^2 t_k^2} dp$ derived above for the imaginary part of the corresponding integral expressing the contribution of a Riemann zero. Each term in the r.h.s. of (3) is definitely positive which proves the

THEOREM 5. An equality

$$\int_0^\infty \frac{\arg(\varsigma(1/2+\alpha+it)(-1/2+\alpha+it))}{(1/2-\alpha)^2+t^2} dt = -\int_{1/2+\alpha}^\infty \frac{\ln(\varsigma(\sigma)(\sigma-1))}{(\sigma-1)(\sigma-2\alpha)} d\sigma,$$ where $\alpha$ is an arbitrary real number $0 \le \alpha < 1/2$, holds true for some $\alpha$ if and only if there

are no Riemann function zeroes with $\sigma > 1/2 + \alpha$. For $\alpha = 0$ this equality is equivalent to the Riemann hypothesis.

REMARK 2. Results of Wang [2] are a particular case of our approach obtained studying the contour integral $I = \int_C \frac{\ln(\varsigma(z))}{a^2 - (z-b)^2} dz$ with $a=0$. Evidently, for this case the pole $t=0$ is located on the contour and the appearing here integral $-i \int \frac{\ln(\varsigma(b+it))}{t^2} dt$ diverges in the vicinity of $t=0$. By this reason Wang was obliged to modify a contour adding a semicircle which excludes the point $z=b$ from it (see Fig. 2). Each zero of the Riemann function lying in the interior of the contour contributes $-2\pi i \int_0^{\sigma_k - b} \frac{dp}{-(p + it_k)^2} = \frac{-2\pi i}{\sigma_k - b + it_k} + \frac{2\pi}{t_k}$ to the contour integral value. The sum $\sum_{\sigma_k > b} \frac{2\pi}{t_k} = 0$ due to the symmetry properties of the Riemann function while in the sum $\sum_{\sigma_k > b} \frac{-2\pi i}{\sigma_k - b + it_k}$ we can pair-wise combine the terms with $it_k$ and $-it_k$ thus obtaining $\sum_{\sigma_k > b, t_k > 0} \frac{-4\pi i (\sigma_k - b)}{(\sigma_k - b)^2 + t_k^2}$.

This is Wang's result. Note also that for not to complicate his somewhat cumbersome expressions further, Wang selected the radius $R$ of a semicircle large enough (actually he used $b=1/2$ and $R=1$) to remove the simple pole of the Riemann function at $z=1$ from the interior of the contour, see Fig. 2. Otherwise, if $b + R < 1$, it is necessary to add $2\pi i \int_R^{1-b} \frac{dp}{-p^2} = 2\pi i \left( \frac{1}{1-b} - \frac{1}{R} \right)$, which is the contribution of this pole, to the integral value.

## 3.3 Integrals with $g(z) = \dfrac{1}{(a^2 - (z-b)^2)^{3/2}}$

It is also instructive to consider $g(z) = \dfrac{1}{(a^2 - (z-b)^2)^{3/2}}$. Due to the reasons discussed above we will limit ourselves with the particular case $g(z) = \dfrac{1}{((1/2-\alpha)^2 - (z-1/2-\alpha)^2)^{3/2}}$ (as before, more general case does not pose problems), where $\alpha$ is an arbitrary real number $0 < \alpha < 1/2$, put $f(z) = \ln(\varsigma(z)(z-1))$ and consider our usual contour $C$ formed by the vertices $b - iY$, $b + iY$, $b + X + iY$, $b + X - iY$ with $1 > b > 1/2$, $X, Y \to +\infty$. Now the function $g(z)$ is double-valued and hence to apply the same reasoning as above we modify our contour indenting the point $z=1$ thus creating the contour $C''$ (see Fig. 3a). An integral over an infinitesimally small circle describing the point $z=1$ vanishes (on this circle the function under the integral sign is asymptotically $Const / \sqrt{z-1}$), the value of the contour integral taken along the "external" three straight lines containing at least one of the points $b + X + iY$, $b + X - iY$ again vanishes when $X, Y \to \infty$, while the function $F(z)g(z)$ changes the sign after indenting the point $z=1$. Thus if RH holds true we have the following identity:

$$-i\int_{-\infty}^{\infty} \frac{\ln(\varsigma(1/2 + \alpha + it) \cdot (-1/2 + \alpha + it))}{((1/2 - \alpha)^2 + t^2)^{3/2}} dt - 2\int_{1}^{\infty} \frac{\ln(\varsigma(\sigma) \cdot (\sigma - 1))}{(1-\sigma)^{3/2}(\sigma - 2\alpha)^{3/2}} d\sigma = 0.$$

Both terms occurring here are purely imaginary and hence

$$\int_{0}^{\infty} \frac{\ln|\varsigma(1/2 + \alpha + it) \cdot (-1/2 + \alpha + it)|}{((1/2 - \alpha)^2 + t^2)^{3/2}} dt = -\int_{1}^{\infty} \frac{\ln(\varsigma(\sigma) \cdot (\sigma - 1))}{(\sigma - 1)^{3/2}(\sigma - 2\alpha)^{3/2}} d\sigma \qquad (5).$$

REMARK 3. It might be helpful to use the example #4.293.14 from Gradshtein and Ryzhik book [8]:

$$\int_0^\infty \frac{x^{\mu-1}\ln(\delta+x)}{(\delta+x)^\nu}dx = \delta^{\mu-\nu}B(\mu,\nu-\mu)[\psi(\nu)-\psi(\nu-\mu)+\ln\delta],$$ that is

$$\int_0^\infty \frac{\ln(|-1/2+\alpha+it|)}{((1/2-\alpha)^2+t^2)^{3/2}}dt = \frac{1}{2(1/2-\alpha)^2}[\psi(3/2)+\gamma+2\ln(1/2-\alpha)],$$ to remove

$\ln(|-1/2+\alpha+it|)$ from the integral sign in the l.h.s. of (5). Thus we obtain

$$\int_0^\infty \frac{\ln|\varsigma(1/2+\alpha+it)|}{((1/2-\alpha)^2+t^2)^{3/2}}dt =$$

$$-\int_1^\infty \frac{\ln(\varsigma(\sigma)\cdot(\sigma-1))}{(\sigma-1)^{3/2}(\sigma-2\alpha)^{3/2}}d\sigma - \frac{1}{2(1/2-\alpha)^2}[\psi(3/2)+\gamma+2\ln(1/2-\alpha)].$$ Here, as usual, $\psi(3/2)=\frac{d}{dx}\ln\Gamma(x)|_{x=3/2}$ is a logarithmical derivative of the gamma-function and $\gamma$ is the Euler constant. In view of $\psi(3/2)=\psi(1/2)+2=2-2\ln 2-\gamma$, we have

$$\int_0^\infty \frac{\ln|\varsigma(1/2+\alpha+it)|}{((1/2-\alpha)^2+t^2)^{3/2}}dt =$$

$$-\int_1^\infty \frac{\ln(\varsigma(\sigma)\cdot(\sigma-1))}{(\sigma-1)^{3/2}(\sigma-2\alpha)^{3/2}}d\sigma - \frac{1}{2(1/2-\alpha)^2}[2-2\ln 2+2\ln(1/2-\alpha)]$$

Introducing now the "semi-contour" $C'''$ (see Fig. 3b) which "semi-indents" the point $z=1$, we can repeat the approach given in the previous section and consider the real part of corresponding integral equality to obtain that if the RH holds true, then:

$$\int_0^\infty \frac{\arg(\varsigma(1/2+\alpha+it)(-1/2+\alpha+it))}{((1/2-\alpha)^2+t^2)^{3/2}}dt + \int_{1/2+\alpha}^1 \frac{\ln(\varsigma(\sigma)\cdot(\sigma-1))}{(1-\sigma)^{3/2}(\sigma-2\alpha)^{3/2}}d\sigma = 0 \quad (6).$$

The second term here is an integral over the finite range which apparently can be calculated very precisely what, we believe, makes this equality especially interesting.

As follows immediately from our theorem 1, to give explicit expressions generalizing (5) and (6) without assuming the truth of the RH, one needs to find the values of the integrals $\int_{b+it_k}^{\sigma_k+it_k} \frac{dz}{(a^2-(z-b)^2)^{3/2}}$, where $a=1/2-\alpha$, $b=1/2+\alpha$ and then to sum up them for all Riemann function zeroes with $\sigma_k > b$. Evaluation of these integrals can be performed in the same manner as an analysis of the integrals $\sum_{\sigma_k>b} \int_{b+it_k}^{\sigma_k+it_k} \frac{dz}{a^2-(z-b)^2}$ given above. Using the same variable change $z = p+b+it_k$ we obtain that each term in this sum is equal to $\int_0^{\sigma_k-b} \frac{\exp(3i\varphi/2)}{(a^2-p^2+t_k^2+4p^2t_k^2)^{3/4}}dp$ where $\varphi = \arctan\left(\frac{2pt_k}{a^2-p^2+t_k^2}\right)$. The conditions $t_k \gg p$, $0 \leq p$ again guaranty that the real part of the above integral is always positive as well as that its imaginary part is positive for $t_k > 0$. This proves the

THEOREM 6. An equality

$$\int_0^\infty \frac{\ln|\varsigma(1/2+\alpha+it)\cdot(-1/2+\alpha+it)|}{((1/2-\alpha)^2+t^2)^{3/2}}dt = -\int_1^\infty \frac{\ln(\varsigma(\sigma)\cdot(\sigma-1))}{(\sigma-1)^{3/2}(\sigma-2\alpha)^{3/2}}d\sigma,$$ where $\alpha$ is an arbitrary real number $0 \leq \alpha < 1/2$, holds true for some $\alpha$ if and only if there are no Riemann function zeroes with $\sigma > 1/2+\alpha$. For $\alpha = 0$ this equality is equivalent to the Riemann hypothesis.

and

THEOREM 7. An equality

$$\int_0^\infty \frac{\arg(\varsigma(1/2+\alpha+it)(-1/2+\alpha+it))}{((1/2-\alpha)^2+t^2)^{3/2}}dt + \int_{1/2+\alpha}^1 \frac{\ln(\varsigma(\sigma)\cdot(\sigma-1))}{(1-\sigma)^{3/2}(\sigma-2\alpha)^{3/2}}d\sigma = 0,$$ where $\alpha$ is an arbitrary real number $0 \le \alpha < 1/2$, holds true for some $\alpha$ if and only if there are no Riemann function zeroes with $\sigma > 1/2+\alpha$. For $\alpha = 0$ this equality is equivalent to the Riemann hypothesis.

*3.4. The case b<1/2: equalities involving infinite sums over the Riemann zeroes on the critical line*

We have already noted that all equalities and expressions given above hold also for $b=1/2$, that is for $\alpha = 0$. Exactly integrals with $b=1/2$ were considered earlier by Wang [2] and Balazard – Saias – Yor [4] (their equality is $\int_{-\infty}^\infty \frac{\ln(|\varsigma(1/2+it)|)}{1/4+t^2}dt = 0$, which is our Theorem 4 for $\alpha = 0$ supplemented with $\int_{-\infty}^\infty \frac{\ln(|1/2+it|)}{1/4+t^2}dt = 0$).

More interesting cases can be obtained when putting $b<1/2$. For this case an integration contour contains the critical line $z=1/2+it$ and hence the sums over infinite sets of the Riemann zeroes do appear. This immediately pose the question of the convergence of the corresponding sums, but actually their convergence is guaranteed by the known properties of the Riemann zeroes. All above results can be generalized for the case $b</1/2$, but here we present only the case

$$\frac{a}{\pi}\int_{-\infty}^\infty \frac{\ln(|\varsigma(1/2-\alpha+it)|)}{a^2+t^2}dt = \ln\left|\frac{\varsigma(1/2-\alpha+a)\cdot(a-1/2-\alpha)}{a+1/2+\alpha}\right| \\ + \sum_{\sigma_k>b}\ln\left|\frac{a+\sigma_k-1/2+\alpha+it_k}{a-\sigma_k+1/2-\alpha-it_k}\right| \quad (7)$$

which is a generalization of (2), and its modification

$$\frac{a}{\pi}\int_{-\infty}^{\infty}\frac{\ln(|\varsigma(1/2-\alpha+it)\cdot(-1/2-\alpha+it)|)}{a^2+t^2}dt = \ln(|\varsigma(1/2-\alpha+a)\cdot(a-1/2-\alpha)|)+$$

$$\sum_{\sigma_k>b}\ln\left|\left(\frac{a+\sigma_k-1/2+\alpha+it_k}{a-\sigma_k+1/2-\alpha-it_k}\right)\right|$$

As usual, $0<\alpha<1/2$. In particular, if $a=1/2-\alpha$ we have

$$\frac{(1/2-\alpha)}{\pi}\int_{-\infty}^{\infty}\frac{\ln(|\varsigma(1/2-\alpha+it)\cdot(-1/2-\alpha+it)|)}{(1/2-\alpha)^2+t^2}dt = \ln(|\varsigma(1-2\alpha)\cdot 2\alpha|)+$$

$$\sum_{\sigma_k>b}\ln\left|\left(\frac{\sigma_k+it_k}{1-\alpha-\sigma_k-it_k}\right)\right| \qquad (8).$$

It is instructive to see how this result can be obtained by a modification of Balazard-Saias-Yor approach [4]. They proved their formula (see above) using the complex variable transform $z=\frac{s-1}{s}$, which transforms to the interior of the circle $|z|<1$ the half plane $\operatorname{Re}s>1/2$, and then applying the generalization of Jensen theorem (which holds true for the Riemann function because the latter is a function belonging to an appropriate Banach space of analytical functions, see [4, 9]) to the circle $|z|=1$.

Instead, the variable transform $z=\frac{s-(1-2\alpha)}{s}$ can be used. Then $|z|=1$ is realized on the straight line $s=1/2-\alpha+it$ and the corresponding transformation transforms to the interior of the circle $|z|<1$ the half plane $\operatorname{Re}s>1/2-\alpha$. Note that $z=0$ corresponds to $s=1-2\alpha$, hence $\ln|f(0)|=\ln|\varsigma(1-2\alpha)\cdot 2\alpha|\neq 0$. Now repeating what is done in [4] (selecting $f(z)=\varsigma(s)(s-1)$ and applying the generalized Jensen theorem) and noting that $\frac{1}{|z|}=\frac{|s|}{|s-(1-2\alpha)|}$ hence $dz=-\frac{1-2\alpha}{s^2}ds$, we have:

$$\int_{-\infty}^{\infty}\frac{\ln|\varsigma(1/2-\alpha+it)(1/2+\alpha-it)|}{(1/2-\alpha)^2+t^2} = \frac{2\pi}{1-2\alpha}(\ln|\varsigma(1-2\alpha)\cdot 2\alpha|+\sum_{\rho}\ln\frac{|\rho|}{|\rho-(1-2\alpha)|})$$

This is nothing else than our relation (8) which is a particular case of (7).

## 4. Numerical results

To support our theoretical considerations we perform some numerical experiments. Particularly, we test (2) for two cases. In the first case we set *a=1* and *b=3/4*, so there is no known zeroes lying inside the contour *C* and we can omit the sum over the Riemann zeroes. In the second experiment (*a=1* and *b=1/4*) we also have to consider the contributions of the zeroes. In both cases we execute the calculations until *t=300* (corresponding to the first 137 zeroes) and until *t=1000* (i.e. including the first 649 zeroes). We perform our calculations using the function NIntegrate of the mathematical software "Mathematica" and the Riemann zeroes tabulated by Odlyzko [10]. Table 1 summarizes the results obtained from the numerical experiments. In Fig. 4 we also give a plot of the integrand appearing in (2). We note the highly oscillatory behaviour of this function.

|  | *t=300* | *t=1000* |
| --- | --- | --- |
| Case 1: *a=1, b=3/4* | l.h.s. = 0.16330508251202144<br>r.h.s. = 0.16330187363718995<br>$\Delta$ = 3.208874831489572×10$^{-6}$ | l.h.s. = 0.16330204796176637<br>r.h.s. = 0.16330187363718995<br>$\Delta$ = 1.7432457641297638×10$^{-7}$ |
| Case 2 *a=1, b=1/4* | l.h.s. = -0.4004296118960703<br>r.h.s. = -0.4004296115785665<br>$\Delta$ = -3.175038010283515×10$^{-10}$ | l.h.s. = -0.39881726259737776<br>r.h.s. = -0.3988170532634759<br>$\Delta$ = -2.0933390187938272×10$^{-7}$ |

Table 1: Numerical results

## REFERENCES

[1] A. Ivic, The Riemann zeta-function, John Wiles & sons, New York, 1985.

[2] F. T. Wang, A note on the Riemann Zeta-function, Bull. Amer. Math. Soc. 52 (1946), 319.

[3] V. V. Volchkov, On an equality equivalent to the Riemann hypothesis, Ukrainian Math. J., 47, 1995, 422.

[4] Balazard M., Saias E. and Yor M. (1999), Notes sur la fonction de Riemann, Advances in Mathematics, 143, 284.

[5] D. Merlini, The Riemann magneton of the primes, Chaos and complexity Lett., Nova Sci. Publ., New York, Vol. 2., N. 1, 2006, 93.

[6] S. Beltraminelli, D. Merlini, and S. K. Sekatskii, A hidden symmetry related to the Riemann hypothesis with the primes into the critical strip, arXiv, math.NT, 0803.1508, 2008.

[7] E. C. Titchmarsh, The theory of functions, Oxford, Oxford Univ. Press, 1939.

[8] I. S. Gradshtein et I. M. Ryzhik, Tables of integrals, series and products, Academic, New York, 1990.

[9] K. Hoffman, Banach spaces of analytical functions, New Yoork, Dover, 1988.

[10] A. Odlyzko, http://dtc.umn.edu/~odlyzko/zeta_tables/zeros1

S. K. Sekatskii, Laboratoire de Physique de la Matière Vivante, IPMC, BSP 408, Ecole Polytechnique Fédérale de Lausanne, CH1015 Lausanne-Dorigny, Switzerland.

E-mail : serguei.sekatski@epfl.ch

S. Beltraminelli, CERFIM, Research Center for Mathematics and Physics, PO Box 1132, 6600 Locarno, Switzerland.

E-mail: Stefano.beltraminelli@ti.ch

D. Merlini, CERFIM, Research Center for Mathematics and Physics, PO Box 1132, 6600 Locarno, Switzerland.

E-mail: merlini@cerfim.ch

**Figure captions**

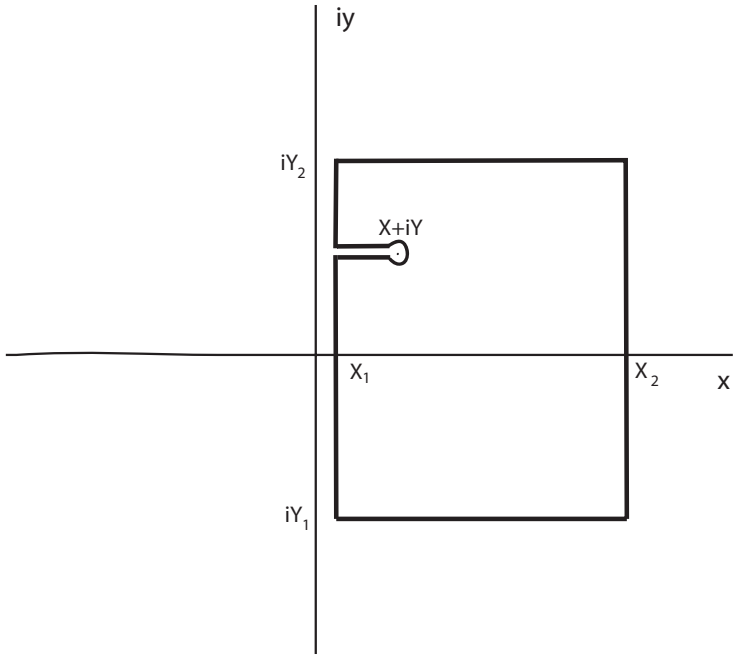

Fig. 1. Contour *C'* using during the proof of the generalized Littlewood theorem.

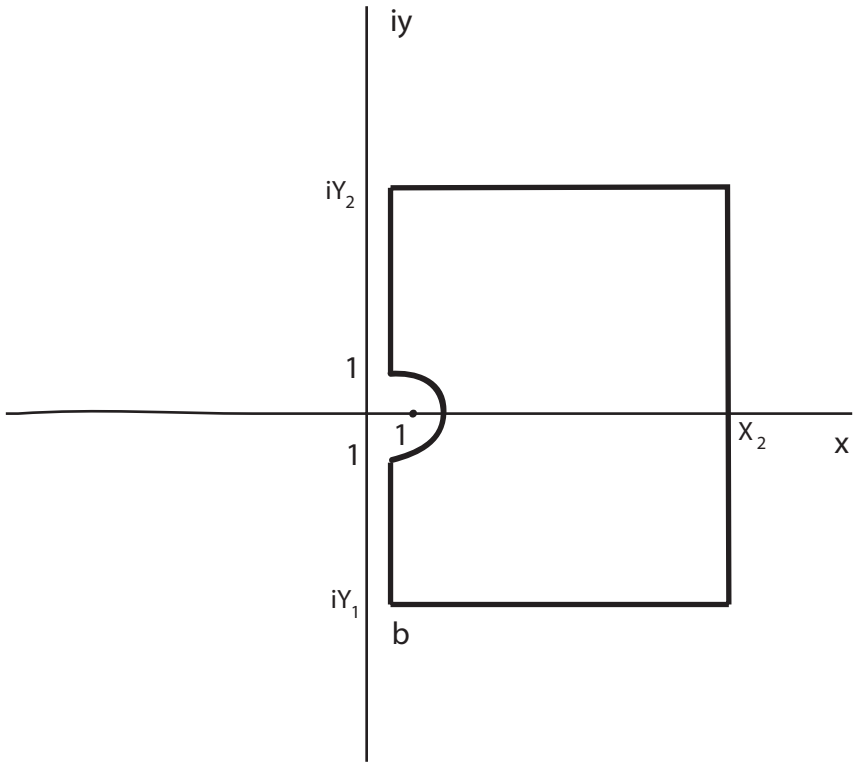

Fig. 2. Contour used during the proof of Wang theorem [2].

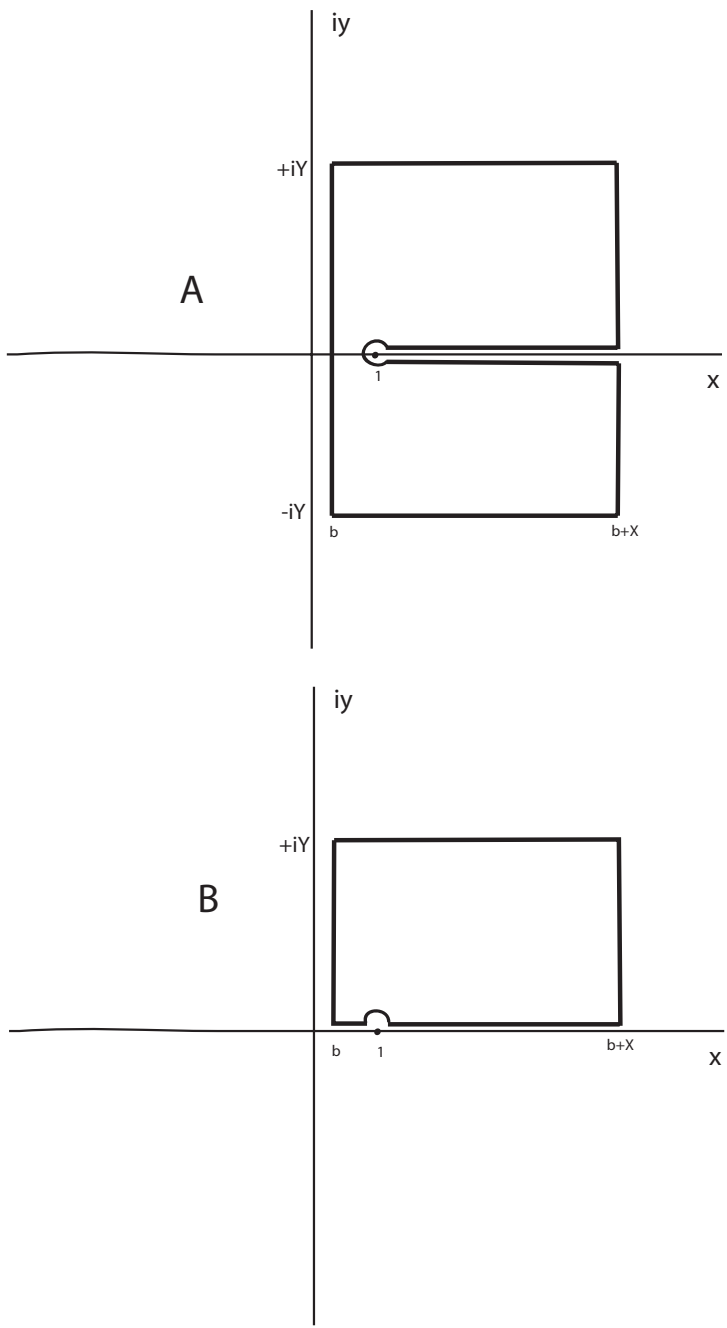

Fig 3. A. Contour $C''$ using for investigation of the integrals $\int_{C''} \frac{\ln(\varsigma(z)\cdot(z-1))}{((1/2-\alpha)^2 - (z-1/2-\alpha)^2)^{3/2}}$. B. "Semi-contour" $C'''$ using for investigation of the same integral. Here $b = 1/2 + \alpha$, $0 < \alpha < 1/2$.

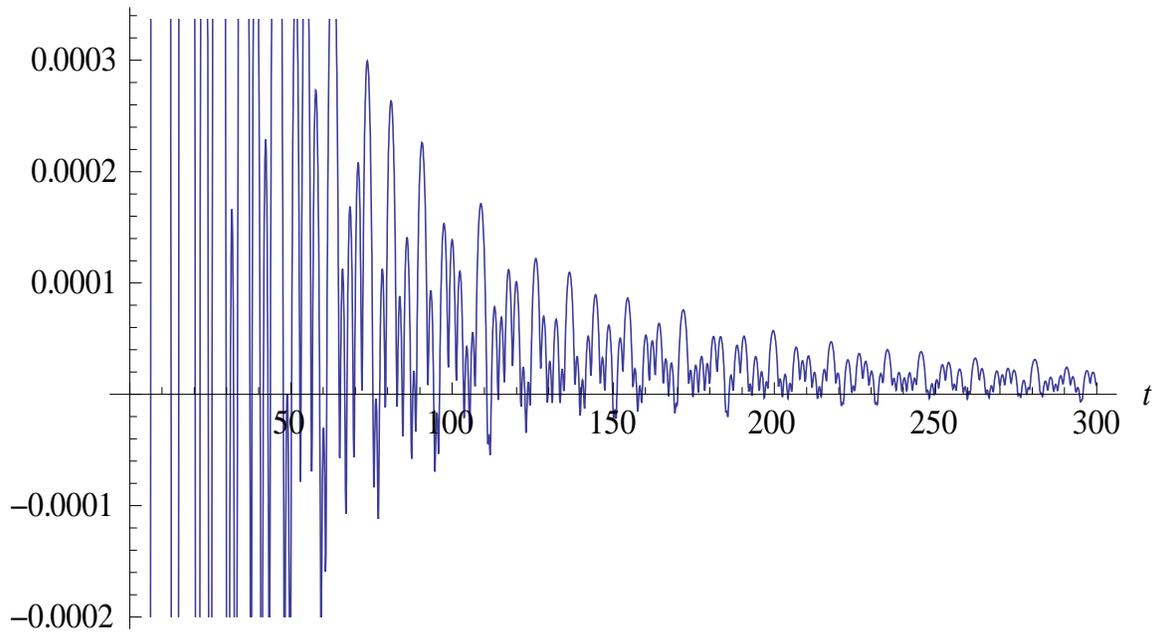

Fig. 4. The integrand of the formula (2) of the paper for *a=1* and *b=1/4*.